\documentclass[12pt]{article}
\usepackage{amsfonts}
\usepackage{xcolor}

\usepackage[plainpages,backref,urlcolor=blue]{hyperref}
\usepackage{tikz}

\title {\bf{Geometry of the Sasakura bundle}}

\author{\bf{Cristian Anghel}}  
\newtheorem{defn}{Definition}[section]

\newtheorem{tth}[defn]{Theorem}
\newtheorem{rem}[defn]{Remark}
\newtheorem{prop}[defn]{Proposition}
\newtheorem{cor}[defn]{Corollary}

\newtheorem{que}[defn]{Question}
\newtheorem{conj1}[defn]{Hartshorne Conjecture}

\begin{document}

\maketitle

\bigskip
\noindent
{\small {{\bf ABSTRACT:}  The Sasakura bundle is a relatively recent appearance in the world of remarkable vector bundles on projective spaces. 
 In fact, it is connected with some surfaces in $\mathbb P^4$ which  missed in early classification papers. 
  The aim of this review is to present various, scattered in the literature, aspects concerning the geometry of this bundle. 
  The last part will be devoted to the place of this bundle in the classification of globally generated locally free sheaves with $c_1 \leq 4$ on $\mathbb P^n$  in a joint paper with I. Coanda and N. Manolache. } }

\bigskip
\noindent
2010 \textit{Mathematics Subject Classification}: 14J60, 14H50.

\

\noindent
Keywords: projective space, vector bundle, elliptic conic bundle.

\tableofcontents

\section{Introduction}

The construction of interesting vector bundles on the projective space is a problem with a long history. It's difficulty is due partially to the great homogeneity of $\mathbb P^n$ in contrast with the case of other varieties, for examples of fibered or product type. The building blocks of general bundles with respect to the operation of direct sum are of course the indecomposable ones:

\begin{defn}
A vector bundle is indecomposable if it is not a direct sum of bundles of smaller ranks.
\end{defn} 
\noindent
It is difficult to construct indecomposable algebraic bundles on $\mathbb P^n$ of small rank $r<n$ for $4\leq n$. In fact the Hartshorne Conjecture asserts that such bundles did not exists for large dimension $n$ and small rank $r$. The precise form in the rank $2$ case is the following:

\begin{conj1}
On $\mathbb P^n$ with $7\leq n$ every rank $2$ bundle splits as a direct sum of line bundles.
\end{conj1} 

\noindent
Some famous examples are:

 -the $r=2$ Horrocks-Mumford bundle on $\mathbb P^4$ $('73)$.
 
 -the $r=n-1$ Vetter-Trautmann-Tango bundle on $\mathbb P^n$ $('73-'76)$.
 
 -the $r=n-1$ nullcorrelation bundle on $\mathbb P^n$ for odd $n$ $('77)$. 
 
 -the $r=3$ Horrocks bundle on  $\mathbb P^5$ $('78)$.
 
 -the $r=3$ Sasakura bundle on $\mathbb P^4$ $('86)$.\\
\noindent
More recent examples are the weighted Tango bundles in arbitrary dimension introduced by Cascini \cite{cas} in ('01), and the bundles of 
Kumar-Peterson-Rao \cite{kpr} in ('02) for low dimensions and various characteristics. The aim of this review is to explain various constructions of the above mentioned Sasakura bundle. The plan of the paper goes as follows: in the next two sections we will describe the original construction of the Sasakura bundle, using \cite{ads} and \cite{sa}. In particular, we will see it's connection with an interesting elliptic conic bundle in $\mathbb P^4$ and with the Beilinson's monads. In the fourth section we will present a beautiful purely geometric construction of the elliptic conic bundle due to Ranestad \cite{ra}. The fifth section is devoted to a totally different viewpoint due to Kumar, Petersen and Rao from \cite{kpr}. In the last section we will describe the Sasakura bundle from the perspective of a joint project with Coanda and Manolache \cite{acm}, aiming to classify globally generated bundles with small first Chern class on $\mathbb P^n$.

\section{The Abo-Decker construction}

The first construction we discuss is due to Abo-Decker-Sasakura in a paper on arxiv:alg-geom/9708023 which appeared one year later \cite{ads}. 
 The aim of the paper, apart the construction of the Sasakura bundle itself, was the construction of an elliptic conic bundle in $\mathbb P^4$:
\begin{defn}
An elliptic conic bundle is smooth projective surface in $\mathbb P^4$ with a projection on a smooth elliptic curve, such that the general fiber is a smooth conic.
\end{defn}
The interest for 
such a surface was due to the following result \cite{ok}: 

\begin{tth}[Okonek]
A smooth surface $X\subset \mathbb P^4$ of degree $8$, sectional genus $5$ and irregularity $1$ is an elliptic conic bundle with exactly $8$ singular fibers composed by pairs of $-1$ lines.
\end{tth}
In fact Okonek claimed that such surfaces does not exists.\\
The main tool in the Abo-Decker-Sasakura construction is the Beilinson theorem which describes coherent sheaves as the cohomology of a certain complex. For a sheaf $\cal F$ on $\mathbb P^4$ one denote by $F^i:=\oplus H^{i+j}({\cal F}(-j))\otimes \Omega ^j(j)$, the direct sum being over all $j$'s. The following theorem from \cite{bei} is a fundamental 
tool in the study of vector bundles on $\mathbb P^n$:
\begin{tth}[Beilinson]
The $F^i$'s forms an increasing complex, exact except in dimension $0$, where the cohomology is $\cal F$.
\end{tth}
\noindent
A first step in the construction is the calculus of the cohomology table for ${\cal J}_X $, the ideal of ( the apriori  hypothetical elliptic conic bundle) $X$:

$ \ \ \ \ \ \ \ \ \ \ \ \ \ \ \ \ \ \ \ \ \ \ \ \ \ \ \ i \uparrow $

$ \ \ \ \ \ \ \ \ \ \ \ \ \ \ \ \ \ \ \ \ \ \ \ \ \ \ \ 4  \mid $

$ \ \ \ \ \ \ \ \ \ \ \ \ \ \ \ \ \ \ \ \ \ \ \ \ \ \ \ \ \ \mid 1 \ 1$

$ \ \ \ \ \ \ \ \ \ \ \ \ \ \ \ \ \ \ \ \ \ \ \ \ \ \ \ \ \ \mid \ \ \ \ 1 \ 1$

$ \ \ \ \ \ \ \ \ \ \ \ \ \ \ \ \ \ \ \ \ \ \ \ \ \ \ \ \ \ \mid \ \ \ \ \ \ \ \ \ 6$

$ \ \ \ \ \ \ \ \ \ \ \ \ \ \ \ \ \ \ \ \ \ \ \ \ \ \ \ \ \ -------- \longrightarrow m$

\noindent
where in the $(m,i)$ box is depicted $h^i({\cal J}_X(m))$. In particular, \\ 
$Ext^1({\cal J}_X(3), {\cal O}(-1))$ is $4$-dimensional, and if it is denoted by $W$, then the identity in $W^*\otimes W$ defines an extension $\cal G$, which is locally free by a generalized  Serre correspondence:
 $$0\rightarrow 4{\cal O}(-1)\rightarrow {\cal G} \rightarrow {\cal J}_X(3) \rightarrow 0.$$ It is the rank $5$ version of the Sasakura bundle.
\noindent
For ${\cal G}(-3)$, the above construction show the following cohomology table:

$ \ \ \ \ \ \ \ \ \ \ \ \ \ \ \ \ \ \ \ \ \ \ \ \ \ \ \ i \uparrow $

$ \ \ \ \ \ \ \ \ \ \ \ \ \ \ \ \ \ \ \ \ \ \ \ \ \ \ \ \ \ \mid $

$ \ \ \ \ \ \ \ \ \ \ \ \ \ \ \ \ \ \ \ \ \ \ \ \ \ \ \ \ \ \mid 1 \ 1$

$ \ \ \ \ \ \ \ \ \ \ \ \ \ \ \ \ \ \ \ \ \ \ \ \ \ \ \ \ \ \mid \ \ \ \ 1 \ 1$

$ \ \ \ \ \ \ \ \ \ \ \ \ \ \ \ \ \ \ \ \ \ \ \ \ \ \ \ \ \ \mid \ \ \ \ \ \ \ \ \ $

$ \ \ \ \ \ \ \ \ \ \ \ \ \ \ \ \ \ \ \ \ \ \ \ \ \ \ \ \ \ ----- \longrightarrow m$ 

\noindent
By Beilinson theorem, $\cal G$ is the cohomology of the monad:
$$0\rightarrow \Omega ^3 (3)\rightarrow \Omega ^2 (2)\oplus \Omega ^1 (1) \rightarrow {\cal O} \rightarrow 0.$$
Moreover, if $\alpha, \beta $ are the maps in the monad, and $e_0,...,e_4$ is a basis in $V$-the underlying vector space of $\mathbb P^4$, using the identification $$Hom(\Omega ^i(i),\Omega ^j(j))\simeq \Lambda ^{i-j}V$$ it can be proved that:
$$\alpha = \left(\begin{array}{c} e_4 \\ e_0 \wedge e_2+e_1 \wedge e_3 \end{array}\right)$$
$$\beta = ( e_0 \wedge e_2+e_1 \wedge e_3  -e_4 ).$$
\noindent
As a conclusion of the above computations, the authors of \cite{ads} obtains the following result:
\begin{tth}
An elliptic conic bundle in $\mathbb P^4$ determines an unique, up to isomorphism and linear change of coordinates, rank-$5$ bundle $\cal G$ given by the monad above.
\end{tth}
\noindent
Moreover, starting with such a $\cal G$, then ${\cal G}(1)$ is globally generated and by a result of Banica \cite{ban} 
the dependency locus of $4$ general sections is a smooth surface with the desired geometric properties.

\section{The Sasakura method}

A different idea was used more than 10 years earlier by Sasakura in $'86$ \\ to construct a stable rank $3$ bundle $\cal E$.
 It is related to the $\cal G$ above by:
\begin{center}
$0\rightarrow 2{\cal O}\rightarrow {\cal G}(1)\rightarrow {\cal E} \rightarrow 0$.
\end{center}
 Also, it is globally generated and the dependency locus of two general sections is an elliptic conic bundle.\\ 
The general setting of the construction is the following:
 $\cal E$ is a rank $r$ vector bundle on $\mathbb P^n$ with first Chern class $c_1$. $s_1,...,s_l$ and $\sigma _1,...,\sigma _k$ are generators of ${H^0}_*(\cal E)$ and ${H^0}_*(\cal E^\vee )$, where the lower star means the direct sum of the respective cohomology of all the twists of the bundle in question:
$${H^i}_*({\cal E}) = \bigoplus H^i({\cal E}(j))$$
 for all $j \in \mathbb Z$.\\
 The generators $s_1,...,s_l$ and $\sigma _1,...,\sigma _k$ above, determines morphisms, the first one surjective  and the second injective
$$ \cal L \rightarrow \cal E \rightarrow  \cal K ,$$ where $\cal L $ and $\cal K $ are direct sums of line bundles.
 In particular, the composition $S: \cal L \rightarrow \cal K $ is a matrix of homogenous polynomials. \\
The construction works in the hypothesis that $\cal K $ has a direct summand of the form $r{\cal O}(m)$. Let $\sigma _{i_1},...,\sigma _{i_r}$ the corresponding generators. 
 By projecting on $r{\cal O}(m)$, we obtain an analogous sequence as above,
$$ {\cal L} \rightarrow {\cal E} \rightarrow r{\cal O}(m),$$ the composition being denoted by $T$.  
 Therefore $\cal E$ is a sub-sheaf in $r{\cal O}(m)$, and moreover a sub-bundle outside the divisor of the form $f:=\sigma _{i_1} \wedge ...\wedge  \sigma _{i_r}$ of degree $r \cdot m-c_1$. 
 Conversely, we can start with a pair $(T,f)$ and ask for conditions under which the resulting sheaf $\cal E$ is 
locally free. Let $I$ the ideal defined by the maximal minors of $T$. The following condition is sufficient for the local freeness of $\cal E$:
the ideal $(I:f)$ define the empty set in $\mathbb P^n$.\\
An useful choice is $f=f'^{r-1}$ with $f'=x_1.....x_{c_1}$ and $T=(T',T'')$ with $T'=f'\cdot Id$ and a convenient $T''$. Using this idea, in  \cite{ads} and \cite{sa} are produced many known bundles, eg. nullcorelation bundle on $\mathbb P^3$, the Horrocks-Mumford on $\mathbb P^4$ and also the Sasakura rank $3$ bundle on $\mathbb P^4$.
 This last one, is constructed with $f=f'=x_0...x_4$ and a convenient {but complicated} $(3\times 8)$ matrix $T$ of forms of degree $4$.

\section{The Ranestad construction}

In $('99)$, Ranestad \cite{ra} gives a purely geometric construction for the elliptic conic bundle.
In fact, his construction works only for elliptic conic bundles (recall they are not minimal) with an 
elliptic quintic scroll in $\mathbb P^4$ as minimal model:
\begin{defn}
An elliptic quintic scroll in $\mathbb P^4$ is a ruled surface over an elliptic curve, embedded with degree $5$ and all fibers as lines.
\end{defn}
\noindent
The main steps are as follows: \\
{\bf Step 1} For every elliptic quintic scroll, denoted $X_5$, in $\mathbb P^4$ there are smooth curves $G$ and $L$ such that $G$ is a rational normal curve intersecting the scroll in $10$ points, and $L$ is a secant of $G$ and a fiber of the scroll.
The ideal of $G\cup L$ is generated by $5$ quadrics, which by an ancient result of Semple $('29)$ defines a Cremona transformation $\varphi $ on $\mathbb P^4$. \\
{\bf Step 2} The restriction of $\varphi $ to $X_5$ is defined by the system $\mid 2H-L\mid $ with $8$ base points: those where $G$ meet the scroll and are not on $L$. \\
{\bf Step 3} This restrction is an embedding of $\hat X_5$, the blow-up at the $8$ points as soon as the $10$ points in $G\cap X_5$ are distinct and no other secant of $G$ is a fiber of the scroll. The image $\varphi (\hat X_5)$ will be of course an elliptic conic bundle.

\

\noindent
Conversely, for a generic elliptic conic bundle, the minimal model is an elliptic quintic scroll and a first intricate construction, produce a quadratic surface $X_2$.
Secondly, a cubic scroll $X_3$ is produced using the secant variety of $G$, the curve above associated with the scroll.
Finally, the converse Cremona is defined by the cubic hypersurfaces through $X_2 \cup X_3$ and it sends the conic bundle to the scroll.  

\begin{rem}
Cf. \cite{ra} one should note that the above constructions works only for elliptic conic bundles with an elliptic quintic scroll as minimal model AND generic with this condition. It is unknown if every such conic bundle can be obtained from its quintic scroll by this construction. Also, it is unknown if every elliptic conic bundle has a quintic scroll as minimal model. 
\end{rem}

\section{The Kumar-Peterson-Rao construction}

In \cite{kpr} the above authors gives new methods for the construction of low rank vector bundles on projective spaces $\mathbb P^4$ and $\mathbb P^5$ in various characteristics. Also, their method was generalized in arbitrary dimension in \cite{bah}. Among some new examples they rediscovered all the known remarkable vector bundles {\bf except} the rank $3$ Sasakura bundle.
 The general construction starts with a sort of Maruyama's elementary transformation:
let $X$ be a projective variety (later it will be $\mathbb P^4$) and $Y$ the divisor of a section $s$ in ${\cal O}_X(Y)$ (later it will be the thickening of order $t$ of a hyperplane in $\mathbb P^4$):
\begin{defn}
If $H$ is a hyperplane in $\mathbb P^4$ with equation $h=0$, the standard thickening of order $t$ of $H$ is the sub-scheme defined by the equation $h^t=0$.
\end{defn}
Let's consider on $Y$ an exact sequence
$$0\rightarrow A \rightarrow F \rightarrow B \rightarrow 0$$ of vector bundles, such that $F$ extends to an $\cal F$ on $X$.
 Let $\cal G$ the kernel of the induced surjection ${\cal F} \rightarrow B$. 

\noindent The next ingredients are two bundles ${\cal L}_1$ and ${\cal L}_2$ on the ambient $X$, 
their restrictions to $Y$, $L_1$ and $L_2$ together with a surjection $L_1 \rightarrow A$ and an injection as vector bundle $B \rightarrow L_2$.
 Suppose that the induced $\phi : L_1 \rightarrow F$ and $\psi : F \rightarrow L_2$ also extends to $\Phi : {\cal L}_1 \rightarrow \cal F$ and $\Psi : {\cal F} \rightarrow {\cal L}_2$.
 Then $\Psi \Phi$ vanishes on $Y$ and one can construct a map
\begin{center}
$\Delta : {\cal F}(-Y)\oplus {\cal L}_1 \rightarrow {\cal F}\oplus {\cal L}_2(-Y)$
\end{center}
given by the matrix below:
$$\alpha = \left(\begin{array}{cc} s\cdot I & \Phi \\ \Psi & s^{-1}\cdot \Psi \Phi \end{array}\right)$$
where $I$ is the identity of $\cal F$ (and $s$ the section which determines $Y$). 

\noindent The role of the above map $\Delta $ whose image is in fact $\cal G$ is tied with the fact that if $\cal F$ splits, under some aditional hypothesis, this splitting will produce sub or quotient bundles of $\cal G$, and consequently lower rank bundles.   
 For example, concerning the sub-bundle case, the authors have the following
\begin{prop}
Suppose:\\
1. $\cal F $ split as ${\cal N}\oplus {\cal N}'$ with induced splitting $F=N\oplus N'$\\
2. there is $\theta :N(-Y) \rightarrow A$ with lift $\Theta :{\cal N}(-Y) \rightarrow {\cal L}_1$ such that $$\Phi \Theta :{\cal N}(-Y) \rightarrow \cal F$$ has image in ${\cal N}'$.\\
Then, there is an induced map $${\cal N}(-Y) \rightarrow \cal G$$ which is a bundle inclusion iff the restriction $$N(-Y)\rightarrow {\cal G}\mid _Y$$ is.
\end{prop}
Roughly speaking, we need decomposable bundles on $Y$ which extends with the decomposition on the ambient $X$.\\
 This goal is achieved through the so called four generated rank two bundles. 
\begin{defn} A vector bundle $B$ on a scheme $Y$ is four generated if there is a completely split (i.e. direct sum of line bundles) rank $4$ bundle $F$ and a surjection $F\rightarrow B$. \end{defn}
 The result below constructs plenty (but of course not all) of four generated rank $2$ bundles on $\mathbb P^3$:
\begin{prop}
Let $T, U, V, W$ a regular sequence of forms of positive degrees $t, u, v, w$ such that $t+w=u+v$ and $r\geq 2$ an integer. Then there is an exact sequence 
$$L_1 \rightarrow F \rightarrow L_2,$$ with maps $\phi $ and $\psi $ such that:\\
1. the bundles above are completely split,\\
2. the images $A$ and $B$ of $\phi $ and $\psi $ are four generated rank $2$ bundles on $\mathbb P^3$.
\end{prop}

\noindent Finally, one can put together these ideas  in the case of $\mathbb P^4$:\\
- one starts with the above $F$ and four generated rank $2$ bundles $A, B$ on $\mathbb P^3$\\
- one consider a $t$-thickening $Y$ of $\mathbb P^3 \subset \mathbb P^4$ and one takes the pull-backs of $F, A$ and $B$ on $Y$\\
- one apply the Maruyama type construction obtaining a rank $4$ bundle $\cal G$ on $\mathbb P^4$\\ 
- for convenient values of the parameters ( the forms $T, U, V, W $ and integers $t, r \geq 2$ ) the bundle $\cal G$ has line sub or quotient bundles\\
- by taking the quotient or the kernel one arrive at rank $3$ bundles on $\mathbb P^4$.

\

\noindent As we can observe from above, this method is quite general and we end this section by the following  natural question:
\begin{que}
Can the Kumar-Peterson-Rao method or its extension from \cite{bah} be used to construct the Sasakura bundle?
\end{que}

\section{The Coanda-Manolache method}

An interesting problem studied in the lasts years is the description of globally generated bundles in projective spaces. This last section aims to describe the place of the Sasakura bundle in this context, following the joint paper with Coanda and Manolache \cite{acm}. Chronologically, the main contributions in 
this area are: \\
 - Chiodera and Ellia  \cite{chel} in ('12) determined the Chern classes of rank $2$ globally generated vector bundles with $c_1\leq 5$ on $\mathbb P^n$,\\
 - Sierra and Ugaglia \cite{siug1} in ('09) classified globally generated vector bundles with $c_1\leq 2$ on $\mathbb P^n$,\\
 - Sierra and Ugaglia  \cite{siug2} in ('14) and independently Manolache with the author \cite{am} in ('13) classified globally generated vector bundles with $c_1\leq 3$ on $\mathbb P^n$,\\
 - Coanda, Manolache and the author \cite{acm} ('13) classified globally generated vector bundles with $c_1\leq 4$ on $\mathbb P^n$.

\

\noindent For example, in the  $c_1\leq 3$ case, the main result from the joint paper with Manolache \cite{am}, can be formulated as follows:
\begin{tth}
Let $E$ an indecomposable globally generated vector bundle on $\mathbb P^n$, with $n\geq 2$, $1\leq c_1\leq 3$ and $H^i(E^*)=0$ for $i=0,1$. Then one of the following holds:\\
- $E={\cal O}(a)$\\
- $E=P({\cal O}(a))$\\
- $n=3$ and $E=\Omega (2)$\\
- $n=4$ and $E=\Omega (2)$\\
- $n=4$ and $E=\Omega^2 (3)$
\end{tth}
where the $P$-operation above for a globally generated bundle $E$, means the dual of the kernel of the evaluation map
$$ H^0(E)\otimes {\cal O} \rightarrow E.$$ 

\noindent The technical condition $H^i(E^*)=0$ for $i=0,1$ is  irrelevant but very useful. In fact any globally generated bundle can be obtained by one which verify this condition by taking the quotient with a trivial sub-bundle and then adding a trivial summand:
\begin{prop}
For any $E$ there is an $F$ which satisfy $H^i(F^*)=0$ for $i=0,1$ such that if $t=h^0(E^*)$ and $s=h^1(E^*)$ then $E\simeq {F \slash s\cal O}\oplus t\cal O$.
\end{prop}
\noindent Another important observation is that for a globally generated bundle one has $c_2 \leq {c_1}^2$. This shows via the standard sequence $$0 \rightarrow (r-1){\cal O} \rightarrow E \rightarrow {\cal J}_Y(c_1)$$ that globally generated bundles with $c_1\leq 3$ are related with sub-varieties of low degree (at most $9$) in $\mathbb P^n$, fact which is also crucial in the classification.

\noindent The main result in the joint paper with Coanda and Manolache \cite{acm} is more complicated. In the $c_1=4$ case there are $16$ indecomposable bundles, the last one being the Sasakura's rank $5$ bundle once twisted ${\cal G}(1)$. 
 However one can formulate the following consequence:
\begin{cor}
Let $E$ an indecomposable globally generated vector bundle on $\mathbb P^n$, with $n\geq 4$, $c_1=4$, $r\geq 2$ and $H^i(E^*)=0$ for $i=0,1$. Then $E$ is:\\
-$P({\cal O}(4))$\\
-$\Omega (2)$ or $\Omega^3(4)$ on $\mathbb P^5$\\
-${\cal G}(1)$.
\end{cor}
where again, the $P$-operation above means the dual of the kernel of the evaluation map. 

\

\noindent The main steps in the classification for $c_1=4$ are of course much more complicated than in the $c_1=3$ case, but can be summarized as follows: \\
 - first we classify the bundles on $\mathbb P^2$ (easy) and $\mathbb P^3$ (hard),\\
 - next we try to decide which bundles can be extended to higher dimensional projective spaces using Horrocks method of killing cohomology. 
 The problem is tied to the theory of varieties of small degree as we have seen, and also to the theory of rank $2$ reflexive sheaves on $\mathbb P^3$ via the exact sequence 

$$0 \rightarrow (r-2){\cal O} \rightarrow E \rightarrow {\cal E}' \rightarrow 0.$$
The twisted rank-$5$ Sasakura bundle ${\cal G}(1)$ appear for $n=4$ and $c_2=8$ and it has the following description:
 consider the surjection $$4{\cal O}(-1)\oplus {\cal O}(-2) \rightarrow {\cal O}$$ defined by $x_0,...,x_3,{x_4}^2$, and the Koszul complex $(C_p, \delta _p)$ for it. Denote by $E'$ the co-kernel of 
$$ \delta _4(4) :   {\cal O}\oplus 4{\cal O}(-1)\rightarrow 4{\cal O}(1)\oplus 6{\cal O}.$$
Then $E={\cal G}(1)$ is the kernel of a surjection $E'\rightarrow {\cal O}(2)$ such that $$H^0(E'(-1))\rightarrow H^0({\cal O}(1))$$ is injective.

\noindent As a last remark, the usefulness of the Koszul complex is not singular in the above case; it appear also in another case, $n=3$ $c_2=8$ where one consider the complex associated with $x_0, x_1, {x_2}^2, {x_3}^2$ and the $E$ is the cohomology of the monad
$${\cal O}(-1) \rightarrow 2{\cal O}(2)\oplus 2{\cal O}(1)\oplus 4{\cal O} \rightarrow {\cal O}(3).$$
Also, other different types of Koszul complexes were used by Kumar, Peterson and Rao to produce interesting deformations of known bundles \cite{kpr1}.

\bigskip

\noindent
Cristian Anghel \\
Department of Mathematics\\
Institute of Mathematics of the Romanian Academy\\
Calea Grivitei nr. 21 Bucuresti Romania\\
email:\textit{Cristian.Anghel@imar.ro}

\end{document}